\documentclass[12pt]{article}
\usepackage{amsmath, amsthm, amscd, amsfonts, amssymb, graphicx, color,cite}
\usepackage{mathrsfs}
\usepackage[T1]{fontenc}
\usepackage[T2C]{fontenc}
\usepackage[T2A]{fontenc}
\usepackage[bookmarksnumbered, colorlinks, plainpages]{hyperref}
\newtheorem{theorem}{Theorem}[section]

\newtheorem{proposition}[theorem]{Proposition}

\theoremstyle{definition}
\newtheorem{definition}[theorem]{Definition}

\theoremstyle{remark}
\newtheorem{remark}{Remark}[section]
\numberwithin{equation}{section}

\begin{document}
	
	\setcounter{page}{1}
	
	
	\begin{center}
		{\LARGE \textbf{A study of Geodesic $(E,F)$-preinvex Functions on Riemannian Manifolds}}
		
		\bigskip

		{\large {Ehtesham Akhter$^{a}$} and {Musavvir Ali$^{b,*}$}}\\
		\textbf{}  \textbf{}\\

		{\small $^{a,~b}$  Department of Mathematics,\\ Aligarh Muslim University, Aligarh-202002, India}
		
		\noindent
		{
			E-mail addresses: $^{b,*}$\url{musavvirali.maths@amu.ac.in}\footnote{${b}$: Corresponding Author},\\
			$^{a}$ \url{ehteshamakhter111@gmail.com}}
	\end{center}
	\bigskip

	{\abstract
		\noindent
			In this manuscript, we define $(E, F)$-invex set, $(E,F)$-invex functions, and $(E,F)$-preinvex functions on Euclidean space. We explore the concepts on the Riemannian manifold. We also detail the fundamental properties of $(E,F)$-preinvex functions and some examples that illustrate the concepts well. We have established a relation between $(E,F)$-invex and $(E,F)$-preinvex functions on the Riemannian manifolds. We introduce the conditions $\mathcal{A}$ and define $(E,F)$-proximal sub-gradient. To explore and demonstrate its applicability to optimization problems, $(E,F)$-preinvex is utilized. In the last, we establish the points of extrema of a non-smooth $(E,F)$-preinvex functions on $(E,F)$-invex subset of the Riemannian manifolds by using $(E,F)$-proximal sub-gradient.}\\
	

	\section{Introduction}
		One of the most crucial aspects of mathematical programming is the search for convexity and generalized convexity because these concepts play a vital role in optimization theory. Several notable generalizations of convexity have been made over the past few decades; see \cite{Pini,Kadakal,Alizadeh,Almutairi,Hudzik,Alomari,Shi,Wang,Liao}. Youness gave the concept of $E$-convex sets and $E$-convex functions also he explored the $E$-convex programming problems in his paper \cite{Youness}. X.M. Yang and X. Chen found some results wrong in paper \cite{Youness}, and then they modified those results (see \cite{Chen,X}). Jian \cite{Jian} came up with the concept of $(E,F)$-convex sets and functions with exciting applications with the generalizations. 
	
	Recently, \cite{Saleh W} Wedad Saleh explored the concept of $(E,F)$-convex sets and functions on the Riemannian manifold and developed some fruitful results in his paper. 

One can find the clear idea and development of the geodesic theory of convexity in the books authored by C. Udriste \cite{Udriste} and T. Rapcs{\'a}k \cite{Rapcsak}. And, for more details on the $E$-convex sets and functions, we refer [\citen{E-convex sets,Saleh W}], from where mainly the idea of $(E,F)$-convex sets, and functions was generated. Convexity and its generalizations on the Riemannian manifolds have gotten more attention from researchers during the last few decades, and a vast amount of literature has been developed on the subject (see \cite{Ahmad1,Agarwal,Agarwal1,Barani,Chen1,Khan,Mititelu,Zhou,Ansari,Saleh}).

	Motivated by the above works and results, particularly in \cite{Saleh W,Jian}, we give the definitions of $(E,F)$-invex sets, $(E,F)$-invex functions and $(E,F)$-preinvex functions. We explore this concept on Riemannian manifolds which is the generalization of $(E,F)$-convexity sets and functions. We establish a relation between $(E,F)$-invex and $(E,F)$-preinvex functions on the Riemannian manifolds. We introduce the conditions $\mathcal{A}$ and define $(E,F)$-proximal sub-gradient. To explore and demonstrate its applicability on problems of optimization, $(E,F)$-preinvexity is applied. In the last section of the present paper, we have obtained extremum points for a non-smooth $(E,F)$-preinvex function defined on $(E,F)$-invex subset of the Riemannian manifold, and that become possible due to use of $(E,F)$-proximal sub-gradient.
	\section{Preliminaries}
		\begin{definition}\cite{Saleh W}\label{d1}
		Let $\mathcal{A}$ be a non-empty subset of $\mathcal{R}^n,$ and  $E,F: \mathcal{A} \rightarrow \mathcal{R}^n$ be two maps. For any $ r_1,s_1 \in \mathcal{A}$ and $\mu \in [0,1]$, $F(s_1)+\mu (E(r_1)-F(s_1))\in \mathcal{A}$, then $\mathcal{A}$ is known as an {\it \bf $(E,F)$-convex set.} 
		
			In this section we recall some well-known definition on Riemannian manifolds.\\
			
		\noindent	Let $\mathcal{O}$ be a  smooth Riemannian manifold($C^{\infty}$) it means that $\langle,\rangle_s$, the metric on  the tangent space $T_{s}\mathcal{O}$ of manifold $\mathcal{O}$ induces a norm $\|,\|_s$. If $z$ and $w$ be two points on $\mathcal{O}$ and $\beta:[c,d]\rightarrow \mathcal{O}$ is a piecewise smooth curve connecting $\beta(c)=z$ to $\beta(d)=w$, its length $L(\beta)$ is given by
		$$L(\beta)= \int_{c}^{d}\|\beta^{'}(t)\|_{\beta^{'}(t)}\,dt.$$
		Now, we define
		$$d(z,w)=\inf\{L(\gamma): \mbox{$\beta$ is a piecewise $C^1$ curve connecting $z$ and $w$}\}$$
		for any $z, w \in \mathcal{O}$. The original topology on $\mathcal{O}$ is then induced by a distance $d$. A Levi-Civita connection, denoted by $\nabla_{A}B$, for any vector fields $A,B \in T\mathcal{O}$ also known as a covariant derivative, is known to exist uniquely on every Riemannian manifold. We also recall that a geodesic is a $C^{\infty}$ smooth path $\beta$ whose tangent is parallel along the path $\beta$, i.e., $\beta$ satisfies the equation $\nabla_{\frac{d\beta(t)}{dt}}\frac{d\beta(t)}{dt}=0.$ Any path $\beta$ connecting  $z$ and $w$ in $\mathcal{O}$ such that $L(\beta)=d(z,w)$ is a geodesic, and it is called a minimal geodesic. Recall that for a given curve $\beta:J \rightarrow \mathcal{O}$, a point $t_o \in J$, and a vector $u_o \in T_{\beta(t_o)}\mathcal{O}$, there exists exactly one parallel vector field $U(t)$ along $\beta(t)$ s.t. $U(t_o)=u_o$. Moreover, the mapping defined by $u_o\longmapsto U(t)$ is a linear isometry between the tangent spaces $T_{\beta(t_o)}$ and $T_{\beta(t)}$, for each $t \in J$. We denote this mapping by $P^{t}_{t_o, \beta}$ and we call it the parallel translation from $T_{\beta(t_o)}\mathcal{O}$ to $T_{\beta(t)}\mathcal{O}$ along the curve $\beta$. A simply connected complete Riemannian manifold with non-positive sectional curvature is called a Hadamard manifold.

	\end{definition}
\section{$(E,F)$-invex Sets and $(E,F)$-preinvex Functions}
	\begin{definition}\label{d2}
		Let $E,F: \mathcal{A} \rightarrow \mathcal{R}^n$ be the functions, defined on non-empty subset $\mathcal{A}$ of $\mathcal{R}^n$, and  $ G: \mathcal{R}^n\times \mathcal{R}^n\rightarrow \mathcal{R}^n$ be a bi-function. If $\forall~ r_1,s_1 \in \mathcal{A}$ and $\mu \in [0,1]$,
	$F(s_1)+\mu G(E(r_1),F(s_1))\in \mathcal{A}$, then $\mathcal{A}$ is known as an {\it \bf $(E,F)$-invex set.} 
	\end{definition}
\begin{remark}
	If we take $E=F$ in the definition $(\ref{d2})$, then the set $\mathcal{A}$ becomes an $E$-invex set.
\end{remark}
\begin{remark}
	If we take $G(E(r_1),F(s_1))=E(r_1)-F(s_1)$ in the definition $(\ref{d2})$, then the set $\mathcal{A}$ becomes an $(E,F)$-convex set.
\end{remark}
	\begin{remark}
			If we take $E=F$= identity map in the definition $(\ref{d2})$, then the set $\mathcal{A}$ becomes an invex set.
	\end{remark}
	\begin{remark}
			If we take $E$ \& $F$ be identity map on $\mathcal{A}$ and $G(E(r_1),F(s_1))=E(r_1)-F(s_1)$ in the definition $(\ref{d2})$, then the set $\mathcal{A}$ becomes a convex set.
	\end{remark}
\begin{definition}\label{d3}
	A real-valued function $\mathcal{H}: \mathcal{A} \rightarrow \mathcal{R}$ defined on the $(E,F)$-invex set $\mathcal{A}$, is known as the {\it \bf $(E,F)$-preinvex} function if
	\begin{eqnarray}
		\mathcal{H}(F(s_1)+\mu G(E(r_1),F(s_1)))\leq \mu \mathcal{H}(E(r_1))+(1-\mu)\mathcal{H}(F(s_1))
	\end{eqnarray}
	holds $\forall~ r_1,s_1 \in \mathcal{A}$ and $\forall~ \mu \in [0,1]$.
\end{definition}
\begin{theorem}
		Let $\mathcal{A}$ be a non-empty $(E,F)$-invex subset of $\mathcal{R}^n$. An arbitrary functions  $\mathcal{H}_i,~i=1,2,..,n$ are $(E,F)$-preinvex functions on $(E,F)$-invex set $\mathcal{A}$ and $\alpha_i\geq0, ~i=1,2,...,n$. Then, $\sum_{i=1}^{n}\alpha_i\mathcal{H}_i$ is also $(E,F)$-preinvex function on the $(E,F)$-invex set $\mathcal{A}$.
\end{theorem}
\begin{definition}\label{d4}
	Let $\mathcal{H}$ be a real-valued differentiable function on the open interval $\mathcal{A}$ of $\mathcal{R}^n$. The function $\mathcal{H}$ is known as {\it \bf $(E,F)$-invex} on $A$ if there are the maps $E,F: \mathcal{A} \rightarrow \mathcal{R}^n$ and a bi-function $G: \mathcal{R}^n\times \mathcal{R}^n\rightarrow \mathcal{R}^n$ such that the following inequality holds
	$$\mathcal{H}(E(r_1))-\mathcal{H}(F(s_1))\geq G(E(r_1), F(s_1))^T\nabla \mathcal{H}(F(s_1))$$
	for each $r_1,s_1 \in \mathcal{A}$ and $\mu \in [0,1]$.
\end{definition}
Now, we will explore all the above definitions in Riemannian manifolds.
	\section{Geodesic $(E,F)$-invex Sets and $(E,F)$-preinvex Functions on Riemannian Manifolds}
\begin{definition}\label{d5}
	Let $\mathcal{B}$ be a non-empty subset of the Riemannian manifold $\mathcal{O}$ and $E,F : \mathcal{B} \rightarrow \mathcal{O}, G: \mathcal{O} \times \mathcal{O} \rightarrow T\mathcal{O}$ be functions such that for every $E(r_1), F(s_1) \in \mathcal{O},~ G(E(r_1),F(s_1)) \in T_{F(s_1)}\mathcal{O}$. Then, the subset $\mathcal{B}$ of $\mathcal{O}$ is said to be {\it \bf geodesic $(E,F)$-invex} if for each $E(r_1), F(s_1) \in \mathcal{B}$, there exists exactly one geodesic $\beta_{E(r_1),F(s_1)}: [0,1]\rightarrow \mathcal{O}$ such that
	$$\beta_{E(r_1),F(s_1)}(0)=F(s_1),~ \beta^{'}_{E(r_1),F(s_1)}(0)=G(E(r_1),F(s_1)),~\beta_{E(r_1),F(s_1)}(s) \in \mathcal{B},~ \forall s \in [0,1].$$
\end{definition}
Now, we define the $(E,F)$-invex function on an open geodesic $(E,F)$-invex subset of a Riemannain manifold.
\begin{definition}\label{d6}
	Let $\mathcal{B}$ be a non-empty open subset of the Riemannian manifold $\mathcal{O}$ and $\mathcal{B}$ be a geodesic invex set. A function $\mathcal{H}:\mathcal{B} \rightarrow \mathcal{R}$ is said to be {\it \bf $(E,F)$-invex} on $\mathcal{B}$ if the following inequality holds
	$$\mathcal{H}(E(r_1))-\mathcal{H}(F(s_1))\geq d\mathcal{H}_{F(s_1)}G(E(r_1),F(s_1)), ~~\forall ~E(r_1), F(s_1) \in \mathcal{B}.$$
\end{definition}
\begin{definition}\label{d7}
	Let $\mathcal{B}$ be a non-empty geodesic $(E,F)$-invex subset of the Riemannian manifold $\mathcal{O}$. A function $\mathcal{H}:\mathcal{B} \rightarrow \mathcal{R}$ is said to be $(E,F)$-preinvex  if for each $E(r_1), F(s_1) \in \mathcal{B}$ and $\forall s \in [0,1]$
	\begin{eqnarray}\label{i2}
		\mathcal{H}(\beta_{E(r_1),F(s_1)})(s)\leq s\mathcal{H}(E(r_1))+(1-s)\mathcal{H}(F(s_1))
	\end{eqnarray}
	holds, where $\beta_{E(r_1),F(s_1)}$ is the unique geodesic in the definition \ref{d5}.
	If the above inequality (\ref{i2}) is strict for $E(r_1)\neq F(s_1)$ and $s \in (0,1)$ then, $\mathcal{H}$ is said to be strict $(E,F)$-preinvex function.
\end{definition}

\begin{proposition}
		Let $\mathcal{B}$ be a non-empty geodesic $(E,F)$-invex subset of the Riemannian manifold $\mathcal{O}$. Consider the function $\mathcal{H}:\mathcal{B} \rightarrow \mathcal{R}$, an $(E,F)$-preinvex, then we have
		
		\noindent$(a)$ Every lower level set of $\mathcal{H}$ defined by
		$$P(\mathcal{H},t):=\{E(r_1) \in \mathcal{B} ~|~ \mathcal{H}(E(r_1))\leq t\}$$ is geodesic $(E,F)$-invex set.
		
		\noindent$(b)$ \begin{equation}\label{p1}
			\begin{aligned}
				& {\text{min~$\mathcal{H}$(w)}}\\
				& \text{subject to} ~\text w \in \mathcal{B}
			\end{aligned}
		\end{equation}
	the set $L$ of the solution of the problem $(\ref{p1})$ is a geodesic $(E,F)$-invex set. Also, if $\mathcal{H}$ is strictly $(E,F)$-preinvex function then, $\mathcal{O}$ contains at the most one point.
\end{proposition}
\begin{theorem}
	Let $\mathcal{B}$ be a non-empty open geodesic $(E,F)$-invex subset of the Riemannian manifold $\mathcal{O}$. Consider that $\mathcal{H}:\mathcal{B} \rightarrow \mathcal{O}$ is a differentiable and $(E,F)$-preinvex function. Then, $\mathcal{H}$ is an $(E,F)$-invex function.
\end{theorem}
\begin{definition}\label{c1}
	Let the functions $E,F: \mathcal{O} \rightarrow \mathcal{O}$ be defined on a Riemannian manifold $\mathcal{O}$. Then, we say that the function $G: \mathcal{O} \times \mathcal{O} \rightarrow T\mathcal{O}$ satisfies the condition $\mathcal{A}$, if for every $E(r_1),F(s_1)\in \mathcal{O}$ and for the geodesic $\beta : [0,1] \rightarrow \mathcal{O}$ satisfying $\beta_{E(r_1),F(s_1)}(0)=F(s_1),~ \beta^{'}_{E(r_1),F(s_1)}(0)=G(E(r_1),F(s_1))$, we have\\
	
	\noindent$(\mathcal{A}_1)$ $P^{o}_{s,\beta}[G(F(s_1),\beta(s))]=-sG(E(r_1),F(s_1))$\\
	\noindent$(\mathcal{A}_2)$ $ P^{o}_{s,\beta}[G(E(r_1),\beta(s))]=(1-s)G(E(r_1),F(s_1))$\\
	
	\noindent for all $s \in [0,1]$.
\end{definition}
\begin{theorem}
		Let $\mathcal{B}$ be a non-empty open geodesic $(E,F)$-invex subset of the Riemannian manifold $\mathcal{O}$. Suppose that the function $\mathcal{H}:\mathcal{B}\rightarrow \mathcal{R}$ is differentiable. If $\mathcal{\mathcal{H}}$ is $(E,F)$-invex on $\mathcal{B}$ and $\mathcal{\mathcal{H}}$ satisfying the condition $\mathcal{A}$ then, $\mathcal{\mathcal{H}}$ is $(E,F)$-preinvex function on $\mathcal{B}$.
\end{theorem}
\begin{theorem}
		Let $\mathcal{B}$ be a non-empty open geodesic $(E,F)$-invex subset of the Riemannian manifold $\mathcal{O}$. Suppose the function $\mathcal{H}:\mathcal{B}\rightarrow \mathcal{R}$ is an $(E,F)$-preinvex. If $F(\bar r_1) \in \mathcal{B}$ is a local optimal solution to the problem \begin{equation}\label{p2}
			\begin{aligned}
				& {\text{minimize~$\mathcal{\mathcal{H}}(F(r_1))$}}\\
				& \text{subject to}~ \text F(r_1) \in \mathcal{B},
			\end{aligned}
		\end{equation}
	then $F(\bar r_1)$ is a global minimum in the problem $(\ref{p2}).$  
\end{theorem}
\begin{definition}
	Let the functions $E,F: \mathcal{O}\rightarrow \mathcal{O}$ be defined on Riemannian manifold $\mathcal{O}$, and function $\mathcal{H}: \mathcal{O} \rightarrow (-\infty, \infty]$ be a lower semi-continuous function. An element $\sigma$ in $T_{F(s_1)}\mathcal{O}$ is an $(E,F)$-proximal sub-gradient of $\mathcal{H}$ at $F(s_1)\in \mbox{dom}(\mathcal{H})$, if $\exists$ the numbers $\mu>0$ and $\lambda>0$ s.t.
	$$\mathcal{H}(E(r_1))\geq \mathcal{H}(F(s_1))+\langle \sigma, \exp^{-1}_{F(s_1)}E(r_1))\rangle_{F(s_1)}-\lambda d(E(r_1),F(s_1))^2,$$ ~$\forall~ E(r_1) \in N(F(s_1),\mu),$
	and where $\mbox{dom}(\mathcal{H}):=\{E(r_1) \in \mathcal{O}: \mathcal{H}(E(r_1))< \infty\}$.
	The set of all $(E,F)$-proximal sub-gradient of $\mathcal{H}$ at $F(s_1)\in \mathcal{O}$ is denoted by $\partial_{(E,F)}\mathcal{H}(F(s_1))$ and is called the $(E,F)$-proximal sub-differential of $\mathcal{H}$ at $F(s_1)$.
\end{definition}
\begin{theorem}
		Let $\mathcal{B}$ be a non-empty open geodesic $(E,F)$-invex subset of the Carton-Hadamard manifold $\mathcal{O}$ with $G(E(r_1),F(s_1))\neq0$ for $E(r_1)\neq F(s_1).$ Also, let $\mathcal{H}:\mathcal{B}\rightarrow (-\infty, \infty]$ be a lower semi-continuous $(E,F)$-preinvex function. Assume $F(s_1)\in$ domain$(\mathcal{H})$ and $\sigma \in \partial_{(E,F)}\mathcal{H}(F(s_1)).$ Then, $\exists$ a real number $\mu$ s.t.
		\begin{equation}\label{i5}
			\mathcal{H}(E(r_1))\geq \mathcal{H}(F(s_1))+\langle \sigma,G(E(r_1),F(s_1))\rangle_{F(s_1)},~\forall~E(r_1)\in \mathcal{B}\cap N(F(s_1),\mu ) 
		\end{equation}
\end{theorem}
	

\begin{thebibliography}{33}
			\bibitem{Barani}
		A.	Barani and M.R. Pouryayevali, {\it Invex sets and preinvex functions on Riemannian manifolds}. Journal of Mathematical Analysis and Applications, 328(2), (2007), pp.767-779.
		\bibitem{E-convex sets}
		A. Iqbal, S. Ali and I. Ahmad, {\it On geodesic E-convex sets, geodesic E-convex functions and E-epigraphs}. Journal of Optimization Theory and Applications, 155, (2012), pp.239-251.
			\bibitem{Saleh}
		A.	K{\i}l{\i}{\c{c}}man and W. Saleh, {\it On geodesic strongly E-convex sets and geodesic strongly E-convex functions}. Journal of Inequalities and Applications, 2015(1), pp.1-10.
			\bibitem{Udriste}
		C.	Udriste, {\it Convex functions and optimization methods on Riemannian manifolds} (Vol. 297). Springer Science \& Business Media, 2013.
			\bibitem{Youness}
		E.A. Youness, {\it E-convex sets, E-convex functions, and E-convex programming.} Journal of Optimization Theory and Applications, 102(2), (1999), pp.439-450.
			\bibitem{Shi}
		F.	Shi,  G. Ye, W. Liu and D. Zhao, {\it cr-h-convexity and some inequalities for cr-h-convex function}. Filomat, (2022).
			\bibitem{Wang}
		G.	Wang and Y. He, {\it Generalized convexity of the inverse hyperbolic cosine function}. Miskolc Mathematical Notes, 19(2), (2018), pp.873-881.
			\bibitem{Hudzik}
		H. Hudzik and L. Maligranda, {\it Some remarks on s-convex functions}. Aequationes mathematicae, 48(1), (1994), pp.100-111.
		\bibitem{Ahmad1}
		I. Ahmad, M.A. Khan and A. Ishan A., {\it Generalized geodesic convexity on Riemannian manifolds}. Mathematics, 7(6), (2019), p.547.
			\bibitem{Jian}
		J.B. Jian, {\it On $(E, F)$ generalized convexity}. Int. J. Math. Sci. 2, (2003), pp.121–132.
			\bibitem{Liao}
		J. Liao and T. Du, {\it  On Some Characterizations of Sub-b-s-Convex Functions}. Filomat, 30(14), (2016), pp.3885-3895.
			\bibitem{Zhou}
		L.W. Zhou and N.J. Huang, {\it Roughly geodesic  $ B $-invex and optimization problem on Hadamard manifold}. Taiwanese Journal of Mathematics, 17(3), (2013), pp.833-855.	
		\bibitem{Khan}
		M.A. Khan, I. Ahmad and F.R. Al-Solamy, {\it Geodesic r-preinvex functions on Riemannian manifolds}. Journal of Inequalities and Applications, 2014(1), pp.1-11.
			\bibitem{Alomari}
		M. Alomari,	 M. Darus and S.S. Dragomir, {\it New inequalities of Simpson's type for s-convex functions with applications}. Research report collection, (2009), 12(4).
			\bibitem{Alizadeh}
		M.H. Alizadeh and J. Mak{\'o}, {\it On e-convexity}. Miskolc Mathematical Notes, 23(1), (2022), pp.51-60.
	\bibitem{Kadakal}
M.	Kadakal and {\.I}. {\.I}{\c{s}}can, {\it Exponential type convexity and some related inequalities}. Journal of Inequalities and Applications, 2020(1), pp.1-9.
	\bibitem {Almutairi}
O. Almutairi and A. Kili{\c{c}}man, {\it   Generalized Fej{\'e}r--Hermite--Hadamard type via generalized (h- m)-convexity on fractal sets and applications}. Chaos, Solitons \& Fractals, (2021), 147, p.110938.
	\bibitem{Ansari}
Q.H. Ansari, M. Islam, and J.C. Yao, {\it Nonsmooth convexity and monotonicity in terms of a bifunction on Riemannian manifolds}. Journal of nonlinear and convex analysis, 18(4), (2017), pp.743-762.
\bibitem{Agarwal}
R.P. Agarwal, I. Ahmad, A. Iqbal and S. Ali,  {\it Geodesic G-invex sets and semistrictly geodesic $\eta$-preinvex functions}. Optimization, 61(9) (2012), pp.1169-1174.
	\bibitem{Agarwal1}
	R.P. Agarwal, I. Ahmad, A. Iqbal and S. Ali, {\it Generalized invex sets and preinvex functions on Riemannian manifolds}. Taiwanese Journal of Mathematics, 16(5), (2012), pp.1719-1732.
		\bibitem{Pini}
	R. Pini, {\it Convexity along curves and invexity}. Optimization, 29(4), (1994), pp.301-309.
\bibitem{Chen1}
	S.L. Chen, N.J. Huang and D. O'Regan, {\it Geodesic B-preinvex functions and multiobjective optimization problems on Riemannian manifolds}. Journal of Applied Mathematics, 2014.
	 \bibitem{Mititelu}
	S. Mititelu, {\it Generalized invexity and vector optimization on differential manifolds}. differ. Geom. Dyn. Syst. 3, (2001), pp.21-31.
		\bibitem{Rapcsak}
	T. Rapcs{\'a}k, {\it Smooth nonlinear optimization in $R^n$} (Vol. 19). Springer Science \& Business Media, (2013).
	\bibitem{Chen}
	X. Chen, {\it Some properties of semi-E-convex functions.} Journal of Mathematical Analysis and Applications, 275(1), (2002), pp.251-262.
	\bibitem{X}
		X.M. Yang, {\it On E-convex sets, E-convex functions, and E-convex programming}. Journal of Optimization Theory and Applications, 109(3), (2001), p.699.
	\bibitem{Saleh W}	
W. Saleh,  {\it Hermite–Hadamard type inequality for (E, F)-convex functions and geodesic (E, F)-convex functions.} RAIRO-Operations Research. 56(6), 2022 Nov 1;:4181-9.
	\end{thebibliography}
\end{document}